\documentclass[11pt,reqno,a4paper,twoside]{article}

\usepackage{amsmath,amsthm,amstext,amscd,amssymb,euscript,mathrsfs}
\usepackage{calrsfs}
\usepackage{epsfig}

\newcommand{\R}{\mathbb R}
\newcommand{\N}{\mathbb N}
\newcommand{\bbL}{\mathbb L}

\newcommand{\E}{\mathbb E}
\newcommand{\Zd}{\mathbb Z^d}

\renewcommand{\phi}{\varphi}

\newcommand{\la}{\ensuremath{\Lambda}}
\newcommand{\si}{\ensuremath{\sigma}}

\newcommand{\pee}{\ensuremath{\mathbb{P}}}

\newcommand{\loc}{\mathcal{L}}
\newcommand{\ra}{\rangle}
\newcommand{\lo}{\langle}
\def\1{{\mathchoice {\rm 1\mskip-4mu l} {\rm 1\mskip-4mu l}
{\rm 1\mskip-4.5mu l} {\rm 1\mskip-5mu l}}}

\newtheorem{theorem}{{\small T}{\scriptsize HEOREM}}[section]
\newtheorem{corollary}{{\bf{\small C}{\scriptsize OROLLARY}}}[section]
\newtheorem{proposition}{{\bf{\small P}{\scriptsize ROPOSITION}}}[section]
\newtheorem{lemma}{{\bf{\small L}{\scriptsize EMMA}}}[section]
\newtheorem{remark}{{\bf{\small R}{\scriptsize EMARK}}}[section]
\newtheorem{definition}{{\bf{\small D}{\scriptsize EFINITION}}}[section]

\renewenvironment{proof}[1]
{\noindent{{\bf{\small{ P}{\scriptsize ROOF}}}.}\hspace{0.1cm} #1} {$\;\qed$\newline}

\newcommand{\beq}{\begin{eqnarray}}
\newcommand{\eeq}{\end{eqnarray}}

\newcommand{\ba}{\begin{align*}}
\newcommand{\ea}{\end{align*}}

\newcommand{\be}{\begin{equation}}
\newcommand{\ee}{\end{equation}}

\newcommand{\bl}{\begin{lemma}}
\newcommand{\el}{\end{lemma}}

\newcommand{\br}{\begin{remark}}
\newcommand{\er}{\end{remark}}

\newcommand{\bt}{\begin{theorem}}
\newcommand{\et}{\end{theorem}}

\newcommand{\bd}{\begin{definition}}
\newcommand{\ed}{\end{definition}}

\newcommand{\bp}{\begin{proposition}}
\newcommand{\ep}{\end{proposition}}

\newcommand{\bc}{\begin{corollary}}
\newcommand{\ec}{\end{corollary}}

\newcommand{\bpr}{\begin{proof}}
\newcommand{\epr}{\end{proof}}

\newcommand{\bi}{\begin{itemize}}
\newcommand{\ei}{\end{itemize}}

\newcommand{\ben}{\begin{enumerate}}
\newcommand{\een}{\end{enumerate}}

\newcommand{\caC}{{\mathscr C}}

\newcommand{\caD}{{\EuScript D}}

\newcommand{\caH}{{\mathcal H}}
\newcommand{\caI}{{\mathcal I}}
\newcommand{\caJ}{{\mathcal J}}
\newcommand{\caK}{{\mathcal K}}
\newcommand{\caL}{{\mathcal L}}
\newcommand{\caM}{{\mathcal M}}

\newcommand{\caO}{{\mathcal O}}
\newcommand{\caP}{{\mathcal P}}
\newcommand{\caQ}{{\mathcal Q}}

\newcommand{\caV}{{\mathcal V}}

\newcommand{\tor}{{{\mathbb T}_d^N}}
\newcommand{\avtor}{\frac{1}{|\tor|}\sum_{i\in\tor}}

\newcommand{\dmu}{\dot{\mu}}
\begin{document}
\title{Hamiltonian and Lagrangian for the trajectory of the empirical distribution and the empirical measure of Markov processes}
\author{
Frank Redig$^{\textup{{\tiny(a)}}}$, Feijia Wang$^{\textup{{\tiny(b)}}}$\\
{\small $^{\textup{(a)}}$
Delft Institute of Applied Mathematics,}\\
{\small Technische Universiteit Delft}\\
{\small Mekelweg 4, 2628 CD Delft, Nederland}\\
{\small $^{\textup{(b)}}$ Mathematisch Instituut Universiteit Leiden}\\
{\small Niels Bohrweg 1, 2333 CA Leiden, The Netherlands}}

\maketitle

\begin{abstract}
We compute the Hamiltonian and Lagrangian associated to the large deviations
of the trajectory
of the empirical distribution for independent Markov processes, and of the empirical measure
for translation invariant interacting Markov processes.
We treat both the case of jump processes (continuous-time Markov chains and interacting particle systems)
as well as diffusion processes.
For diffusion processes, the Lagrangian is a quadratic form of the deviation of the trajectory from
the Kolmogorov forward equation. In all cases, the Lagrangian can be interpreted as
a relative entropy or relative entropy density per unit time.
\bigskip

\noindent

\end{abstract}

\section{Introduction}
The Gibbs formalism (DLR equations, variational principle) plays a crucial role in
statistical mechanics of equilibrium systems. Roughly speaking a lattice spin system is called ``Gibbs''
if it can be described in terms of Boltzmann-Gibbs weights with an interaction such that
the total interaction of the spin at the origin with all other spins is finite, uniformly in
all configurations.
Beyond equilibrium, the appearance of
Gibbs measures is less obvious as is illustrated both by the loss of Gibbs property  
in the course of stochastic dynamics of Glauber type shown in \cite{efhr} (modeling heating and cooling),
as well as by the expected non-Gibbsianness of non-equilibrium stationary states.
Outside equilibrium it is natural to think of a Gibbsian description in terms of histories,
i.e., using trajectories of the system. In the context of translation invariant lattice spin systems,
one then ends up naturally with a description of the system modulo translations, i.e., 
on the level of the trajectory of the empirical measure.

In \cite{efhr2} we explained how Gibbs-non-Gibbs transitions in lattice spin systems can be related
to a bifurcation phenomenon for the optimal trajectories -in the sense of large deviations-
of the empirical measure conditioned to arrive at a given measure at a fixed time $T>0$.
In \cite{kuel}, \cite{RW} this idea was developed in the mean-field context, i.e., for
the trajectory of the magnetization.
The idea is to consider a translation invariant stochastic dynamics, and study the trajectory of the empirical measure.
More precisely, one conditions the trajectory to arrive at time
$T>0$ at a given empirical measure $\nu$, and at time zero one
gives a certain cost $i_0(\mu)$ to each translation invariant starting measure $\mu$. This cost is determined by the choice of the Gibbs measure
$\mu^G_0$ from which the dynamics is started, i.e., this cost
equals the relative entropy density $s(\mu|\mu^G_0)$ of $\mu$ w.r.t.\ this initial Gibbs measure.
The total cost of a trajectory arriving at time $T$ at $\nu$ is then the sum of $i_0(\mu)$ with the path space cost
$\Psi(\mu,\nu,T)$ of the optimal (in the sense of large deviations) trajectory starting from $\mu$ at time zero and
arriving at $\nu$ at time $T>0$. The set of minimizers $\caO_T(\nu)$ of this total cost of trajectories arriving at $\nu$ at
time $T$ is then the object which should be investigated
in order to understand whether or not the initial Gibbs measure, evolved over a time $T$ (denoted by $(\mu^G_0)_T$) is a Gibbs measure.

The precise conjecture is the following. If we have uniqueness for every conditioning of the empirical measure at time $T>0$, i.e., if
$\caO_T (\nu)$ is a singleton for
every choice of $\nu$, then this should correspond to Gibbsianness of the time-evolved measure $(\mu^G_0)_T$. Conversely, if we have
non-uniqueness for a particular conditioning of the empirical measure at time $T$ (a so-called ``bad measure''), then we have non-Gibbsianness
of the distribution $(\mu^G_0)_T$ at time $T$.

As stated before, the total cost to arrive at time $T>0$ at a given empirical measure is the sum of the initial cost and a path cost, determined by the Markovian dynamics.
The path cost is usually of the form of a Lagrangian action. This means,
informally written, that the probability of a trajectory of the empirical measure, where one averages
shifts of the point mass of the lattice-spin configuration over the box $[-N,N]^d$,is expected to
behave as
\[
\pee \left((\bbL_N (\si_t))_{0\leq t\leq T}\approx (\mu_t)_{0\leq t\leq T}\right) \approx \exp\left(-(2N+1)^d\int_0^T \Xi (\mu_s,\dot{\mu}_s) \ ds\right)
\]
where $\approx$ has to be interpreted in the sense of the large deviation principle on the space of trajectories of translation invariant probability measures.

The Lagrangian
$\Xi (\mu_s,\dot{\mu}_s)$ is the object we are after in the present paper.

More precisely, we consider two cases in the present paper.

First, in the context of {\em independent Markov processes} on a general state space $E$, we study the Lagrangian of the associated to the large deviations of the trajectory of {\em the empirical
distribution}
\[
\loc_N= \frac1N\sum_{i=1}^N \delta_{X^i_t}
\]
which is a random probability measure on $E$.
We compute explicitly the Hamiltonian and provide information
on the associated Hamiltonian trajectories for finite state space Markov chains.

For diffusion processes, the Lagrangian is a quadratic form associated to the generator.
For Markov chains, the Lagrangian is less explicit (except for two state Markov chains), but can still
be characterized as a relative entropy production.
The study of the large deviations of the trajectory of the empirical distribution has to be considered as {\em the intermediate step} between the magnetization
(studied in \cite{kuel}, \cite{RW}) and the empirical measure. In particular, for finite state space Markov chains,
the empirical distribution is still a finite dimensional object.
The corresponding Gibbs-non-Gibbs transitions associated to uniqueness or non-uniqueness of optimal trajectories are then
situated
in the context of general mean field models, and the notion of Gibbsianness developed there, see e.g. \cite{kuel}.

Second, in the context of {\em translation invariant locally interacting Markov processes}, we consider {\em the trajectory of the empirical
measure}, and compute explicitly the Hamiltonian, both for (interacting) diffusion processes and for jump processes in the class of interacting particle
systems
\cite{ligg}.
In the context of diffusion processes, the Lagrangian is a quadratic form, while in the context of jump processes
(of interacting particle systems type),
the Lagrangian is less explicit, but also in that setting a relative entropy production (density) characterization can be given.

This study is a step in the research programme proposed in \cite{efhr2}. Given the Hamiltonians
and Lagrangians computed in the present paper, one can then characterize bifurcation phenomena,
i.e., non-uniqueness of optimal trajectories for particular choices of initial costs.
We leave this problem for future work and focus here on the explicit form of the Hamiltonian
and Lagrangian. The full and mathematically complete proof of the validity of the trajectory large deviation principle both for the empirical distribution as well as for the empirical measure will be considered in two future works \cite{k}, \cite{kr}. This amounts to prove that
the Hamiltonian and Lagrangian that we compute here correspond to a unique non-linear semigroup, coinciding with the Nisio semigroup associated to the Lagrangian \cite{fk}.

Our paper is organized as follows. In section 2 we give a general computation of the Feng-Kurtz Hamiltonian
for the trajectory of the empirical distribution. In section 3 we study the case of finite continuous-time
Markov chains. In section 4 we consider the case of diffusion processes. In section 5 we consider the
case of interacting Markov processes, both of jump type (interacting particle systems
in the spirit of \cite{ligg}) and of diffusion type.

\section{The trajectory of the empirical distribution: general case}
We consider $\{X_t:0\leq t\leq T\}$ a (Feller) Markov process
on a state space $E$. We assume $E$ to be a locally compact Polish space. Relevant cases for
the present paper are, $E$ a finite set (finite Markov chains), or
$E=\R^k$ or a compact submanifold of $\R^k$ (diffusions, diffusions in a domain). The computation of this section is however
valid for general $E$.

We denote by $Q$ the generator of the process $\{X_t:0\leq t\leq T\}$, i.e.,
\[
Q f(x)= \lim_{t\to 0}\frac1t\left( \E_x f(X_t)- f(x) \right)
\]
for $f\in \caD (Q)$, and where $\E_x$ denotes expectation in the process starting from $X_0=x$.
The corresponding semigroup is denoted by $S_t$, i.e., $S_t f(x)= \E_x f(X_t)$. For $E$ compact
$S_t$ acts on $\caC(E)$, the space of continuous functions, for cases such as
$E=\R^d$, $S_t$ acts on $\caC_0(E)$, the space of continuous functions vanishing at
infinity. We further denote $\caC_b (E)$ the space of bounded continuous
functions on $E$ (of course in the compact case we have $\caC(E)=\caC_b(E)$).
For $\mu$ a finite Borel measure on $E$
and $f\in \caC (E)$, we denote $\lo \mu, f\ra=\int fd\mu$.
We denote by $\caP(E)$ the set of probability measures on $E$.

We now let $\{X^i_t:0\leq t\leq T\}$
be independent copies of the process $\{X_t:0\leq t\leq T\}$
starting at initial points $X^i_0= x_i$, and consider the empirical distribution
\be\label{empdis}
\caM_N(t) =\frac1N\sum_{i=1}^N \delta_{X^i_t}
\ee
This is a random probability measure on $E$, i.e., a random element of $\caP(E)$,
which in the limit $N\to\infty$ converges
to the solution of the Kolmogorov forward equation.

More precisely, if at time zero, $\caM_N(0)\to \mu$ (where $\mu$ is a probability measure on $E$), then
at time $t$, $\caM_N (t)\to \mu_t$, where
$\mu_t$ solves
\be\label{kolmo}
\frac{d\mu_t}{dt}= Q^*\mu_t
\ee
where $Q^*$ denotes the dual generator defined via
\[
\lo \mu, Q f\ra = \lo Q^* \mu, f\ra
\]
for all $f$ in the domain of $Q$.

Indeed, by the law of large numbers, for all $f\in \caC_b(E)$,
\[
\lo \caM_N (t), f\ra=\frac1N\sum_{i=1}^N \E_{x_i}f(X^i_t)\to \int f\ d\mu_t
\]
where $\mu_t= S^*(t)\mu$ denotes the law of $X_t$ when started initially from $X_0$ distributed
according to $\mu$, and where $\to$ denotes convergence almost surely.

The convergence $\caM_N(t)\to\mu_t$ is a manifestation of the law of large numbers, and therefore
it is natural to expect an associated large deviation principle, i.e.,
\be\label{large}
\pee \left(\{\caM_N(t):0\leq t\leq T\}\approx \{\mu_t: 0\leq t\leq T\}\right)
\approx
\exp\left(-N\caI (\{\mu_t: 0\leq t\leq T\})\right)
\ee
Where $\approx$ has to be interpreted in the sense of the large deviation principle \cite{dz}, in
a suitable topology on the space of trajectories, i.e.,
lower bound for open sets $G$ of trajectories
\be\label{largeopen}
\liminf_{N\to\infty}\frac1N\log\pee \left(\{\caM_N(t):0\leq t\leq T\}\in G\right)
\geq -\inf_{\gamma\in G}\caI (\gamma)
\ee
and for upper bound for closed
sets $F$ of trajectories
\be\label{largeclosed}
\limsup_{N\to\infty}\frac1N\log\pee \left(\{\caM_N(t):0\leq t\leq T\}\in F\right)
\leq -\inf_{\gamma\in F}\caI (\gamma)
\ee

By the Markov property, the rate function $\caI$ has the form of a Lagrangian ``action''
\be
\caI \left(\{\mu_t: 0\leq t\leq T\}\right)=\int_0^T \caL (\mu_s, \dot{\mu}_s) \ ds
\ee
where $\dot{\mu}_s$ denotes the weak derivative of the trajectory at time $s$, defined via
\be\label{rororo}
\lo \dot{\mu}_s, f\ra = \frac{d}{ds} \lo \mu_s, f\ra
\ee
Notice that $\dot{\mu}_s$ is certainly well-defined on functions $f$ in \eqref{rororo}
in the domain of the generator $Q$, but can possibly not be extended as a finite signed measure
on the whole space. We leave the formulation of the precise space on which $\dot{\mu}_s$ lives
to the companion papers \cite{k,kr} where full proofs are given.
Our aim here is to compute the Lagrangian $\caL$.

As explained in the introduction, this opens the road to an analysis of bifurcation phenomena related to Gibbs-non-Gibbs
transitions, as is done on the level of the magnetization in \cite{kuel}, \cite{efhr2}, \cite{RW}. The case
of the empirical distribution should correspond to Gibbs-non-Gibbs phenomena in the context of mean-field
models, where the mean field interaction is a function of possibly several empirical averages (rather than only of
the magnetization).

Notice that the {\em expression for the Lagrangian $\caL$} is independent of the
precise topology (on the space of trajectories of probability measures on $E$)
in which the large deviation principle \eqref{large}
holds. Usually, one then first considers the weakest topology which is the
product topology (point-wise convergence at every time), and if one wants to strengthen
the topology to e.g.\ uniform topology, one proves exponential tightness in that
topology.

In this paper we focus on the computation of the lagrangian $\caL$ with
the scheme of Feng and Kurtz \cite{fk}, explained e.g. in \cite{RW}.

In our context this means that we first compute the non-linear generator. To explain this, we need
some more notation. First notice that $(X^1_t,X^2_t, \ldots, X^N_t)$ is a Markov process with
generator
\be\label{indgen}
\caQ_N f(x_1,\ldots, x_n)= \sum_{i=1}^N Q_i f
\ee
where $Q_i$ denotes the generator $Q$ applied to the $i$-th coordinate.

The first computation in the Feng-Kurtz scheme is then the non-linear generator
\be\label{nonlin}
H F (\mu)=\lim_{N\to\infty, \caM_N(x_1,\ldots,x_n)\to\mu} \frac1N\left(e^{-NF(\caM_N)} \caQ_N e^{NF(\caM_N)}\right)
\ee
If $\caH F$ is of the form $\caH (\mu, \nabla F)$, with $\caH$ a strictly convex function
in the second variable, then we call $\caH (\mu, f)$ the Hamiltonian, and
the corresponding Lagrangian is then given by the Legendre transform of $\caH$:
\be\label{lagleg}
\caL (\mu, \alpha) =\sup_{f\in \caC(E)}\left(\int f \ d\alpha - \caH (\mu, f)\right)
\ee
The interpretation of the ``gradient''
$\nabla F$ is straightforward when we are in the context of finite state space Markov chains, because
the set $\caP(E)$ is then finite dimensional. In the context of diffusion processes or more general
Markov processes, the gradient will be a (context dependent) functional derivative.

The second variable of the Lagrangian \eqref{lagleg}
is the velocity variable, which in our context is a signed measure of
total mass zero.

The Hamiltonian $\caH(\mu, f)$ can be obtained as follows:
\be\label{ham}
\caH (\mu,f)=\lim_{N\to\infty, \caM_N\to\mu} \frac1Ne^{-N\lo\caM_N, f\ra} \caQ_N e^{N\lo\caM_N, f\ra}
\ee
Notice here that for a given $f\in \caC (E)$, the function
$e^{N\lo\caM_N, f\ra}=e^{\sum_{i=1}^N f(x_i)}$ is a function from $E^N$ to $\R$, on which
the generator $\caQ_N$ can act, i.e., the notation in
$\caQ_N e^{N\lo\caM_N, f\ra}$ means the generator $\caQ_N$ acting on that
function of $x_1,\ldots,x_N$.

The $\mu$ variable is interpreted as the ``position'' and the $f$ variable as
the ``momentum'' (dual to the velocity variable in the Lagrangian formalism).

By the form \eqref{indgen} of the independent generator, we can compute the Hamiltonian $\caH(\mu, f)$:
\beq\label{hamco}
\caH (\mu,f) &=&\lim_{N\to\infty, \caM_N\to\mu}\frac1N e^{-N\lo\caM_N, f\ra} \caQ_N e^{N\lo\caM_N, f\ra}
\nonumber\\
&=&
\lim_{N\to\infty, \caM_N\to\mu}\frac{1}{N}\sum_{i=1}^N e^{-f(x_i)} Q e^{ f(x_i)}
\nonumber\\
&=&
\int e^{-f} Q e^f \ d\mu
\eeq
Notice that since $\caH (\mu, 0)=0$, for
the corresponding Lagrangian \eqref{lagleg}
we have
\[
\caL (\mu, \alpha) \geq \left(\lo\alpha, 0\ra -\caH (\mu, 0)\right)=0
\]
i.e., the Lagrangian is automatically non-negative (as it should be since it is the integrand of the rate function).
\section{Finite state space continuous-time Markov chains}
\subsection{Hamiltonian and Lagrangian}
In this case $E=\{a_1,\ldots, a_k\}$ is a finite set, of which we denote the elements
by $a,b,\ldots$.
A function $f:E\to\R$ is identified with a column of numbers
$f_a, a\in E$, so we will use both notations $f(a)$, or $f_a$, idem
for probability measures (identified with rows $\mu_a, a\in E$).

The continuous-time Markov chain is defined via its transition rates between states $a,b\in E$,
denoted by $r(a,b)$. The generator is given by
\be\label{fingen}
Q f(a)= \sum_{b\in E} r(a,b) (f(b)-f(a))
\ee
For a probability measure $\mu_a, a\in E$ we then have
the Kolmogorov forward equation for the distribution at time $t$, denoted $\mu_a(t), a\in E$:
\be\label{velo}
\frac{d\mu_a(t)}{dt} = \sum_{b} \left(r(b,a) \mu_b (t)-r(a,b)\mu_a(t)\right)
\ee
with initial condition $\mu_a(0)=\mu_a, a\in E$.

The Hamiltonian \eqref{hamco} is given by
\be\label{finham}
\caH (\mu, f) =\sum_{a,b\in E} \mu_a r(a,b)(e^{f_b-f_a}-1)
\ee
The corresponding Lagrangian is given by
\be\label{finlag}
\caL(\mu, \alpha)= \sup_f \left(\sum_a f_a \alpha_a-\caH(\mu,f)\right)
\ee
The $f=f^*(\alpha)$ which realizes the supremum satisfies
\be\label{optf}
\alpha_b =\sum_{a} \left(\mu_a r(a,b) e^{f^*_b-f^*_a}-\mu_b r(b,a) e^{f^*_a-f^*_b}\right)
\ee
This leads to
\be
L(\mu, \alpha)=\sum_{a,b} \mu_b r(b,a)\left(f^*_ae^{f^*_a- f^*_b}- f^*_b e^{f^*_a-f^*_b}-(e^{f^*_a-f^*_b}-1)\right)
\ee
Defining the ``modified'' rates
\[
r^* (b,a)= r(b,a)e^{f^*_a- f^*_b}
\]
the equation \eqref{optf} reads
\be\label{optf1}
\alpha_b =\sum_{a} \left(\mu_a r^*(a,b)-\mu_b r^*(b,a) \right)
\ee
which can be interpreted as follows. The modified rates are such that
they produce ``velocity'' (\eqref{velo}) equal to $\alpha$, when started from initial measure $\mu$.
In terms of these modified rates $r^*$, the Lagrangian reads
\be
L(\mu, \alpha)=\sum_{a,b} \mu_b r^*(b,a)\log\left(\frac{r^*(b,a)}{r(b,a)}\right)- \sum_{a,b}\mu_b(r^*(b,a)-r(b,a))
\ee
This can be interpreted in terms of relative entropy as follows.
The Radon Nikodym derivative of the path space measure of the process with rates $r^*$ w.r.t.\ the process
with rates $r$ is given by the Girsanov formula
\be\label{fingirs}
\frac{d\pee^{[0,T]}_{r^*}}{d\pee^{[0,T]}_r}=
\exp\left(\sum_{a,b}\left( \log\left(\frac{r^*(b,a)}{r(b,a)}\right) N^{b,a}_T - (r^*(b,a)-r(b,a))T\right)\right)
\ee
where $N^{(b,a)}_T$ denotes the number of transitions from $b$ to $a$ in $[0,T]$
The corresponding relative entropy of these two processes is then given by
\[
s(\pee^{[0,T]}_{r^*}|\pee^{[0,T]}_r)=\int d\pee_{r^*}^{[0,T]} \log\left(\frac{d\pee^{[0,T]}_{r^*}}{d\pee^{[0,T]}_r}\right)
\]
Taking the limit $T\to 0$ in this expression, starting from initial distribution $\mu$, we find
the connection with the Lagrangian:
\be\label{relentlag}
\lim_{T\to 0} \frac1T s(\pee^{[0,T]}_{r^*}|\pee^{[0,T]}_r)=\caL (\mu, \alpha)
\ee

In words this means the following. In order to compute $\caL (\mu, \alpha)$, we have to consider an auxiliary
Markov process with rates that from starting from $\mu$ produce velocity
(in the sense of \eqref{velo}) equal to $\alpha$. The relative entropy of this
process w.r.t.\ the original process in a small interval of time $[0,t]$ is then given by
$t\caL (\mu, \alpha) + O(t^2)$. The Lagrangian $\caL (\mu, \alpha)$ can thus be viewed as
{\em ``relative entropy production'' needed to force the process to have speed $\alpha$ when started from $\mu$}.
In particular for the evolution according to the Kolmogorov forward equation: $\alpha= Q^*\mu$, the cost is of course zero,
and we indeed have in that case $r^*= r$ and $\loc(\mu, Q^*\mu)=0$.

\subsection{Hamiltonian trajectories for finite Markov chains}
The Hamiltonian \eqref{finham} has Hamiltonian trajectories given by
\beq\label{dishameq}
\dot{f}_a &=& -\frac{\partial \caH}{\partial \mu_a}=-\sum_{b} r(a,b) (e^{f_a-f_b}-1)\nonumber\\
\dot{\mu}_a &=& \frac{\partial \caH}{\partial f_a}=\sum_b \left(\mu_b r(b,a) e^{f_a-f_b}-\mu_a r(a,b) e^{f_b-f_a}\right)
\eeq
The interpretation of the second equation is the following.
For a trajectory with ``momentum'' $f$, the motion of the probability measure is that of a Markov process
with rates which are modified according to $f$ via
\be\label{modr}
\tilde{r}(a,b)= r(a,b) e^{f_b-f_a}
\ee
Indeed, for the modified rates $\tilde{r}$, the second equation of \eqref{dishameq} reads simply
\[
\dot{\mu}_a = \sum_{b} \mu_b \tilde{r} (b,a)-\mu_a \tilde{r}(a,b)
\]
which is precisely the Kolomogorov forwards equation for the evolution of a probability distribution in
a Markov chain with rates $\tilde{r}$.

The equation for the momenta, i.e., the first equation of \eqref{dishameq} can be rewritten using
the variables $u_a= e^{f_a}, a\in E$:
\[
\dot{u}_a = -\sum_b r(a,b) (u_b-u_a) = -(Q u)_a
\]
which has the solution
\be\label{momentsol}
u(t)= e^{-tQ} u(0)
\ee

The equation for the ``position variables '' $\mu_a$ is linear and reads
\be\label{poseq}
\mu (t) = M(u(t))\mu (t)
\ee
with $M$ a matrix depending on the solution of the momentum variables, given by
\be\label{posmat}
M_{a,b}(f)= r(b,a) \frac{u_a}{u_b} -\left(\sum_c r(a,c) \frac{u_c}{u_a}\right)\delta_{a,b}
\ee
This matrix has column sums equal to zero, i.e., for all $b\in E$ we have $\sum_{a} M_{a,b}=0$, which corresponds to the conservation
of mass $\sum_{a} \mu_a(t)=1$ in the Hamiltonian evolution. More precisely, the matrix $M_{a,b}$ is precisely
the adjoint of the generator corresponding to the modified rates $\tilde{r}$ defined in
\eqref{modr}.

We thus conclude that the Hamiltonian trajectories are still Markovian, corresponding with
time-dependendent rates, steered by the solution of the momentum equation \eqref{momentsol}.

The solution of \eqref{poseq} is given by
\be\label{optimal}
\mu (t) = e^{\int_0^t M(u(s))} \mu(0)
\ee
which means that we have the form of the optimal trajectories, up to the determination of the integration constants
given by $u(0)$ and $\mu(0)$.
Although the form \eqref{optimal}, \eqref{momentsol}
looks quite explicit, it is not easy in general to find explicit tractable formulas for $\mu(t)$.
The action or path-space cost of an optimal trajectory
\[
\caI =\int_0^T \loc (\mu_s,\dot{\mu}_s)\ ds
\]
can be rewritten in Hamiltonian formalism as
\be\label{hamcost}
\caJ (\{ \mu(s), f(s):0\leq s\leq T\})= \left(\sum_a \int_0^T f_a(t) \dot{\mu}_a(t)\ dt\right) - T\caH(\mu(0), f(0))
\ee
This means that in order to find the optimal cost to travel from a starting measure $\mu(0)=\mu$ towards
a measure $\mu(T)=\nu$ at time $T$, one has to plug in
the solution
\eqref{optimal}, \eqref{poseq} into the expression \eqref{hamcost}, and determine the integration constants
$\mu(0), f(0)$ by initial and final condition.
This leads to a function $\Psi (\mu,\nu,T)$ which is the optimal path cost to travel from $\mu$ to
$\nu$ in time $T$. In concrete situations beyond two state Markov chains, in practice,  this function
is hard to obtain closed formulas for (an issue which we do not want to pursue here).

{\bf Example: Two state symmetric flipping}\\
To see a concrete example of an explicit solution, we consider the case of two states flipping at rate 1, which corresponds
with mean-field independent spin flip
dynamics, treated before in \cite{kuel}, \cite{efhr2}, \cite{RW}.

In that case, the state space is given by $E=\{1,2\}$, the matrix $Q$ is given by
\[
Q=\left(
\begin{array}{cc}
-1 & 1\\
1 & -1
\end{array}
\right)
\]
and the matrix $M$ of \eqref{posmat}
is given by
\[
M=\left(
\begin{array}{cc}
-\frac{u_2}{u_1} & \frac{u_1}{u_2}\\
\frac{u_2}{u_1} & -\frac{u_1}{u_2}
\end{array}
\right)
\]
where $u= (u_1,u_2)^T$ satisfies
\be\label{quu}
\dot{u} = -Q u
\ee
The equation
\[
\dot{\mu} = M\mu
\]
can be differentiated w.r.t.\ time once more, which gives
\[
\frac{d^2\mu}{dt^2}= \left(\frac{dM}{dt} + M^2\right)\mu
\]
Explicit computation, using \eqref{quu}  then gives
\[
\frac{dM}{dt} + M^2=
\left(
\begin{array}{cc}
2 & -2\\
-2 & 2
\end{array}
\right)
\]
which gives the equations
\[
\frac{d^2\mu_1(t)}{dt^2}= 2\mu_1(t)-2\mu_2(t)=-2\frac{d^2\mu_2(t)}{dt^2}
\]
Putting $\mu_1-\mu_2= x$ we have,
\[
\frac{d^2x}{dt^2}= 4x
\]
which gives $x_t = C_1 e^{2t} + C_2 e^{-2t}$ as solutions as found before in
\cite{kuel}, or \cite{efhr2}. From this the optimal trajectory starting at $\mu$
arriving at $\nu$ and its cost can easily be inferred.

\br
The fact that $\frac{dM}{dt} + M^2$ is a constant matrix is quite exceptional. Even in the two state case,
if the rates $r(1,2)=\alpha\not= r(2,1)=\beta$, the matrix $\frac{dM}{dt} + M^2$ is not
constant and differentiating the equation \eqref{poseq} once more does not lead to further simplification.
\er

\section{Diffusion processes}
\subsection{Hamiltonian and Lagrangian}
Here we consider the state space $E=\R^n$ and diffusion processes with generator
\be\label{difgen}
Q= \sum_i b_i(x)\partial_i + \sum_{ij}a_{ij}(x) \partial^2_{ij}
\ee
where
$\partial_i$ denotes partial derivative w.r.t. $x_i$.
Here $b_i(x)$, $a_{ij} (x)$ are supposed to be Lipschitz and sufficiently smooth, ensuring uniqueness of the solution
of the corresponding stochastic differential equation.

The covariance $a_{ij}(x)$ is assumed to be a non-degenerate positive definite matrix, i.e., we assume
that it is bounded from below by a multiple of the identity.

The Hamiltonian $\caH(\mu, f)$ given in \eqref{hamco} can then be computed
and this yields:
\beq\label{difham}
\caH(\mu, f) &=&
e^{-f} Q e^f \ d\mu
\nonumber\\
&=&
\int \left(Qf + \sum_{ij}\partial_i f (x) \partial_j f(x) a_{ij} (x)\right) \ d\mu (x)
\eeq

The measures $\mu$ that we will have to consider are absolutely continuous probability measures
w.r.t.\ Lebesgue measure, $\mu = \mu(x) dx$, where with slight abuse of
notation we use the symbol $\mu$ both for the measure and its density.

Although we are in the infinite dimensional context here, because the Hamiltonian
is a quadratic form, the corresponding Lagrangian can be obtained more easily than in the previous
subsection.

Define the quadratic
form
\be\label{quad}
J_\mu (f,f)= \int\left(\sum_{ij}\partial_i f (x) \partial_j f(x) a_{ij} (x)\right) \ d\mu (x)
\ee
for $f$ in the domain of $A_\mu$. 
\br
Notice that this quadratic form corresponds to the carr\'e du champ operator, i.e.,
\[
J_\mu (f,f)=  \int\Gamma^Q_2 (f,f) d\mu
\]
where
\[
\Gamma^Q_2(f,f)= Qf^2 - 2fQ f
\]
is the carr\'e du champ operator.
\er
To this quadratic form corresponds a positive self-adjoint operator $A_\mu$ (linearly depending on $\mu$) such that
\[
J_\mu (f, f)= \frac12\lo f, A_\mu f\ra
\]
where $\lo f,g\ra=\int f(x) g(x) \ dx$ is the usual $L^2$ inner product.

With this notation, the Hamiltonian can be written in the form
\be
\caH (\mu, f)= \lo \mu, Qf\ra + \frac12\lo f, A_\mu f\ra = \lo Q^*\mu, f\ra + \frac12\lo f, A_\mu f\ra
\ee
Then, the corresponding Lagrangian is computed
\beq\label{genlagdif}
\caL (\mu, \alpha) &=& \sup_{f} \left(\lo f, \alpha\ra- \lo Q^*\mu, f\ra - \frac12\lo f, A_\mu f\ra \right)
\nonumber\\
&=&
\frac12 \lo (\alpha- Q^* \mu), A_\mu^{-1} (\alpha-Q^*\mu) \ra
\eeq

Where $\lo f, A_\mu^{-1} f\ra$ is to be interpreted in the sense of the spectral theorem, i.e.,
$ \|A_\mu^{-1/2}f\|_2^2$
for $f$ in the domain of $A_\mu^{-1/2}$. The Lagrangian is then defined to be infinite
when $(\alpha- Q^* \mu)$ is not in the domain of $A_\mu^{-1/2}$ (cf.\  the abstract form of Schilder's theorem
in abstract Wiener spaces
see \cite{dz}).

We see that the ``typical trajectory'' which follows
the Kolmogorov forward equation has zero cost, since in that case $\dot{\mu}=\alpha= Q^*\mu$, and
hence $\caL(\mu,\alpha)=0$, and in general, the Lagrangian is a quadratic expression in the deviation
of the trajectory
from the Kolmogorov forward equation.

To illustrate this formula, let us consider first the simplest example of the present context, i.e., dimension $n=1$,
drift $b=0$, $a=1/2$, corresponding
to a one-dimensional Brownian motion. The generator is
\[
Q= \frac12 \frac{d^2}{dx^2}
\]
$Q^*=Q$.
The quadratic form \eqref{quad} reads in this case
\[
J_\mu (f,f)=\frac12 \int \mu(x)(f')^2 dx
\]
and the corresponding operator
\[
A_\mu = \frac{d}{dx}\left( \mu(x)\frac{d}{dx}\right)
\]
which gives
\be\label{brownlagra}
\caL (\mu, \alpha) = \frac12 \left\lo\nabla^{-1}\left(\alpha-\frac12 \mu''\right), \frac1\mu\nabla^{-1}\left(\alpha-\frac12 \mu'' \right)\right\ra
\ee
The rigorous  meaning of the formal expression
$\lo\nabla^{-1} f, \nabla^{-1} g\ra$ is the innerproduct in
the space $H_{-1}$, i.e., $\lo (-\Delta)^{-1/2} f, (-\Delta)^{-1/2}g\ra$, with
$\Delta=  \frac{d^2}{dx^2}$, or equivalently the dual space of
the Hilbert space $H_1$.
\br
The rate function \eqref{brownlagra} has also been obtained in the context
of the study of the hydrodynamic limit for independent Brownian particles,
in \cite{ko}. In general, it is an interesting question to understand the relation
between the rate functions which are computed in this paper and the rate functions
for deviations of the hydrodynamic limit, see e.g. \cite{kl}.
For Brownian particles, they coincide because of scale invariance of the Brownian motion.
\er

\subsection{Relative entropy interpretation}

As in the case of finite state space Markov chains, also for diffusion processes, the Lagrangian \eqref{brownlagra} can be interpreted in terms of relative entropy. Let us illustrate this for one-dimensional Brownian motion as the reference process, i.e., $n=1$ with $b=0, a=1/2$
in \eqref{difgen}.

A diffusion process on $\R$ with drift $b(x)$ and variance equal to one has the generator
\[
Q_b= b(x)\frac{d}{dx} + \frac12 \frac{d^2}{dx^2}
\]
if we start this process from a measure $\mu=\mu(x) dx$, then the infinitesimal change
at time zero is given by the adjoint generator working on $\mu$, i.e.,
\be\label{apdri}
\frac12 \frac{d^2\mu (x)}{dx^2} + \frac{d}{dx} (b(x) \mu(x))= (Q_b^* \mu) (x)
\ee
In particular, for $\alpha$, a given absolutely continuous signed measure of total mass zero,
we can find the drift $b$ that corresponds to it by solving the equation
\be\label{ap}
\frac12 \frac{d^2\mu (x)}{dx^2} + \frac{d}{dx} (b(x) \mu(x))=\alpha(x)
\ee

The process with drift $b$ has a corresponding path space measure on the Wiener space of continuous trajectories
denoted by $\pee^{[0,T]}_b$, and we have the Girsanov formula
\be\label{girsdif}
\frac{d\pee_b}{d\pee_0} =\exp \left(\int_0^T b(W_s) dW_s -\frac12 \int_0^T b^2 (W_s) ds\right)
\ee
relating $\pee^{[0,T]}_b$ with the path space measure of the reference process $\pee^{[0,T]}_0$
The relative entropy of the process with drift $b$ w.r.t. the zero drift process is thus given by
\beq
&&s(\pee^{[0,T]}_b|\pee^{[0,T]}_0)= \int d\pee_b \log \left(\frac{d\pee_b}{d\pee_0}\right)
\nonumber\\
&=&
\E^{(0)} \left(\exp \left(\int_0^T b(W_s) dW_s -\frac12 \int_0^T b^2 (W_s) ds\right)\left(\int_0^T b(W_s) dW_s -\frac12 \int_0^T b^2 (W_s) ds\right)\right)
\nonumber\\
\eeq
where the expectation $\E^{(0)}$ is over the standard Brownian motion, i.e., w.r.t. $\pee^{[0,T]}_0$.
Computing then
\[
\lim_{T\to 0} \frac1T s(\pee^{[0,T]}_b|\pee^{[0,T]}_0)
\]
starting from a distribution $\mu$ for the reference process, and using
\[
\E^{(0)}\left(\int_0^T b(W_s) dW_s\right)^2= \E^{(0)}\left(\int_0^T b^2(W_s) ds\right)
\]
gives
\[
\lim_{T\to 0} \frac1T s(\pee^{[0,T]}_b|\pee^{[0,T]}_0)=\frac12\int b^2 (x) \mu(x) dx
\]
which is equal to $\loc(\mu, \alpha)$ given in \eqref{brownlagra},
because by \eqref{ap}
\[
\frac{d}{dx} (b(x) \mu(x))= \alpha -\frac12 \mu''(x)
\]

Hence, as in the finite Markov chain case, we see that the Lagrangian can be interpreted as
the infinitesimal relative entropy cost to produce a ``velocity'' $\alpha$ when started from $\mu$.
In particular, when $\alpha=Q^*\mu$ this cost is zero, corresponding to the fact that
the evolution according to the Kolmogorov forward equation is
an optimal trajectory with zero cost.

\section{Trajectory of the empirical measure}
\subsection{Context and notation}
In the context of translation invariant interacting systems, the empirical distribution is no longer a natural object
because of interactions. In particular, the empirical distribution as a function of time is no longer a Markov process.
The natural object capturing the essential information about the time evolution, modulo translations
is then given by the empirical measure. In order to describe this setting, we need some more notation.
For $N\in\N$ we denote $ V_N = \{ -N,\ldots,N\}^d$ and denote by $\tor$ the $d$-dimensional torus, i.e.,
$V_N$ endowed with addition modulo $2N+1$.

We will consider translation invariant
systems on this torus which for large $N$ have to be thought of as approximations of an infinite interacting system where the individual components
live on the lattice $\Zd$.

The configuration space is $\Omega_N= E^\tor$, where $E$, the single-site space, is a locally compact Polish space.
Further we denote $\Omega= E^{\Zd}$ the state space of the infinite volume process.
As in the previous sections, we mostly consider $E$ or a finite set (interacting particle systems) or $E=\R^n$ (or a submanifold
of $\R^n$) (interacting diffusion processes). Elements of $\Omega_N$ are denoted $\si,\eta,\xi,\ldots$, and for
$\si\in\Omega_N$, $i\in\tor_N$, $\si_i$ denotes the value of the configuration at site $i$.
On $\tor$ we have the addition modulo $2N+1$, and correspondingly, the shift $\tau_i$ defined on $\Omega_N$ via
\be\label{shift}
\left(\tau_i (\si)\right)_j = \si_{j+i}
\ee
on functions
$f:\Omega_N\to\R$ via $\tau_i f(\si)=f(\tau_i\si)$, and on probability measures via $\int f d (\tau_i\mu)=\int \tau_i f d\mu$.
If $A$ is a linear operator on functions $f:\Omega_N\to\R$ then we define its shift over $i$ to be
$\tau_i A \tau_{-i}$, and an operator is called translation invariant if for all $i$, $\tau_i A\tau_{-i}= Q$.
A measure is translation invariant if $\tau_i\mu=\mu$. Natural translation invariant measures on $\Omega_N$ are obtained
by periodizing translation invariant measures on $\Omega$, i.e., starting from $\si$ distributed according to a translation
invariant measure on $\Omega$, we consider
$\si^N_i = \si_i, i\in V_N$, periodically extended to the whole lattice.
Conversely, if we have a probability measure $\mu_N$ on $\Omega_N$ we naturally associate to it
a probability measure on the infinite configuration space $\Omega$, namely we consider the periodic extension
of a configuration drawn from $\mu_N$ to the whole lattice $\Zd$.
This justifies the fact that with slight abuse of notation
we will use sometimes the same symbol $\mu_N$ for a translation invariant measure on $\Omega_N$ as well as for the corresponding
translation invariant measure on $\Omega$.
We denote by $\pee_{inv} (\Omega)$ the set of translation invariant probability measures on $\Omega$.

A function $f:\Omega\to\R$ is called local if it depends  on a finite number of coordinates, i.e., if there
exists a (minimal) finite set $D_f$, called the dependence set of $f$ such that for all $\si,\eta\in\Omega$:
$f(\si_{D_f}\eta_{\Zd\setminus D_f} ) = f(\si)$, i.e., the value of the function is not influenced by changing
the configuration outside $D_f$. Obviously, a local function $f:\Omega\to\R$ can be thought of as being a function
$f:\Omega_N\to\R$ as well, for $N$ large enough such that $V_N\supset D_f$.
The translation $\tau_i f$ of local function is obviously local, with dependence set $D_{\tau_i f} =D_f +i= \{x+i:x\in D_f\}$.

An linear operator (possibly unbounded) $A: \caD(A)\subset \caC(\Omega)\to \caC(\Omega)$ is local
if it acts only on $\eta_i$, for $i$ a finite set $D=D_A\subset\Zd$ of vertices. A local operator
acts naturally on functions $f:\Omega_N\to\R$ for $N$ large enough, such that $D_A$ is contained
in $V_N$.

\subsection{Translation invariant sequence of local generators}
\bd
A translation invariant sequence of local generators is defined to be a
a sequence of generators of the form
$\loc_N = \sum_{i\in \tor} \tau_i Q\tau_{-i}$, with $Q$ a local generator, such that
the corresponding infinite volume generator
$\loc = \sum_{i\in \Zd} \tau_i Q\tau_{-i}$ is well-defined, corresponds to a unique Markov process on $\Omega$,
and has a subset of local functions
as a core. The generator $Q$ is then called the ``source generator''.
\ed
As a consequence, the corresponding processes $\{\si_{N,t}:t\geq 0\}$ converge weakly in path space to
the infinite volume process $\{\si_t:0\leq t\leq T\}$ with generator $\loc$.
Moreover, for the associated semigroups we have that
$S^N_t f\to S_t f$ uniformly as $N\to\infty$, for all local functions $f$.

Let us give some examples in order to make this concept more concrete.
\ben
\item {\bf Independent Markov processes.} For $Q$ a generator of a Markov process on $E$, we define
\[
\loc_N = \sum_{i\in\tor} \tau_i Q_0\tau_{-i}
\]
Where $Q_0$ is the operator $Q$ working on the variable $\si_0$.
Under the process with this generator $\loc_N$ different components evolve independently, as copies of the process
with generator $Q$.
\item {\bf Spin-flip dynamics.} $E$ is finite set (e.g. $E=\{-1,1\}$ for Ising spins), $\theta: E\to E$ a bijection
such that $\theta (a)\not= a$ for all $a\in E$. Furthermore, a local function $r:\Omega\to\R^+$, with dependence set containing
the origin, is given. The local generator is then defined
$Qf (\si) = r(\si) (f(\theta_0 \si) -f(\si))$, where $\theta_0$ means applying $\theta$ to the coordinate $\si_0$ and leaving
all other coordinates unchanged (similarly we denote $\theta_i$).
The corresponding sequence of generators is then given by
\[
\loc_N f(\si) = \sum_{i} ((\tau_i Q \tau_{-i})f)(\si) = \sum_{i\in \tor} r(\tau_i \si) (f(\theta_i\si)- f(\si))
\]
\item {\bf Interacting diffusions.} For $E=\R$ and for a finite set $D\subset\Zd$, we consider the local
generator
\[
Q f(\si) = \left(\sum_{j\in D} \frac{\partial V(\si_D)}{\partial \si_j} \frac{\partial f}{\partial \si_j} \right) + \frac12\frac{\partial^2}{\partial \si_0^2}
\]
and the corresponding
\[
\loc_N =  \sum_{i}\tau_i Q \tau_{-i}f
\]
This represents a system of diffusions, interacting via the potential $V$. E.g. for a nearest neighbor potential $V:\R\to\R$
in $d=1$, the full generator has the form
\[
\sum_i V'(|\si_i-\si_{i-1}|) \left(\frac{\partial}{\partial \si_i}- \frac{\partial}{\partial \si_{i-1}}\right) + \frac12 \frac{\partial^2}{\partial \si_i^2}
\]
corresponding to $D= \{0,1\}$, $V(\si_D) = V(|\si_1-\si_0|)$.
The core for the generator of the infinite volume process is the set of local smooth ($\caC^\infty_0$) test functions.
\item {\bf Local interacting particle systems.} $E$ is a finite
set. For finite
subsets $D_\alpha\subset\Zd$,  a collection of  $T_\alpha: E^{D_\alpha}\to E^{D_\alpha}$
$\alpha\in \{1,\ldots,k\}$, and corresponding rates $c(\alpha, \si)$
we consider the local generator
\[
Q f(\si)= \sum_{\alpha} c(\alpha, \si)( f(T_\alpha\si_{D_\alpha}\si_{D_\alpha^c})-f(\si))
\]
the corresponding local generators then include of course the previous spin-flip case but also
translation invariant spin-exchange (Kawasaki) dynamics, combination of
spin-flip and spin-exchange, etc.
\item {\bf Local averaging.} For $0\in D\subset\Zd$ finite, and $m_D$ a probability
measure on $E^D$, consider
\[
Q f(\si)= r(\si)\int \left(f(\si'_D\si_{D^c})- f(\si)\right) m_D(d\si'_D)
\]
with $r$ a local function.
In words, this means that with rate $r$, the configuration inside
$D$ is replaced by its average over the measure $m_D$.
Examples of this class are the KMP model (a model of heat conduction) \cite{kmp}, or more generally the thermalized BEP process
\cite{cggr}.
\een
\subsection{Trajectory of the empirical measure}
For a configuration $\si\in\Omega_N$, its corresponding empirical measure is defined by
\be\label{empmeas}
\bbL_N (\si) = \frac{1}{|\tor|}\sum_{i\in\tor}  \delta_{\tau_i\si}
\ee
This is a translation invariant probability measure on $\Omega_N$, capturing all information
about $\si$, modulo translations.

For a configuration on the full lattice, $\si\in\Omega$, with a slight abuse of notation
we also denote
\be
\bbL_N (\si) = \avtor \delta_{\tau_i (\si^N)}
\ee
where $\si^N$ is the periodized configuration obtained from $\si$.

If $\mu$ is a probability measure on $\Omega$, which is ergodic under translations, then, by the Birkhoff ergodic theorem, with
$\mu$ probability one
\[
\bbL_N (\si)\to \mu
\]
as $N\to\infty$, and where ``$\to$'' means weak convergence.

If $(\loc_N)_N$ is a translation invariant sequence of local generators, then we have the
associated Markov processes $\si_{N,t}$ with semigroups $S^N_t= e^{t\loc_N}$. For a probability measure
$\mu$ on $\Omega$, let us denote $\mu_t$ to be the distribution at time $t>0$ in the infinite volume process
$\{\si_t:t\geq 0\}$, started at initial state distributed according to $\mu$. By locality of the generator $\loc$, for
$\mu$ ergodic,
we have that $\mu_t$ is ergodic as well and hence
\[
\bbL_N (\si_t)\to \mu_t
\]
weakly, with probability one.
Therefore the random trajectory of translation invariant probability measures
$\{ \bbL_N (\si_t): 0\leq t\leq T\}$ converges, as $N \to\infty$ to the deterministic
trajectory $\{\mu_t: 0\leq t\leq T\}$. This convergence of a random $\caP_{inv} (\Omega)$-valued
trajectory to a deterministic
$\caP_{inv} (\Omega)$-valued trajectory can be thought of as a law of large numbers
(in an infinite dimensional space), and therefore it is natural to ask for an associated large deviation principle.
For spin-flip dynamics, this was studied in \cite{efhr2}. Here we treat the general case of
a translation invariant sequence of local generators.
This will naturally lead to a {\em local non-linear
operator $\caK_Q$ associated to the local source generator $Q$}, which is in the present context of interacting systems
the analogue
of the non-linear operator $e^{-f} Q e^f$ in section 2.

More precisely, we want to identify the ``path space Lagrangian'' (which is in this section is denoted
by $\Xi$) such that,
in $\caP_{inv} (\Omega)$:
\[
\pee\left(\{ \bbL_N (\si_t): 0\leq t\leq T\}\approx \{ \nu_t: 0\leq t\leq T\}\right) \approx
\exp\left({-|\tor| \int_0^T \Xi (\nu_t, \dot{\nu}_t)}\right)
\]
The Lagrangian is now a function
of a translation invariant probability measure and a translation invariant signed measure of total mass zero
(or a more general distribution on the space of functions belonging to the domain of the generator), and as before,
$\approx$ has to be interpreted in the sense of the large deviation principle, in this case, on the space of trajectories with values in the
set $\caP_{inv}(\Omega)$ of translation invariant probability measures on $\Omega$.

\subsection{The Hamiltonian}
In this section  we compute the
Feng-Kurtz Hamiltonian. This Hamiltonian is now a function
from $\caC(\Omega)\times \caP_{inv} (\Omega)$ to $\R$, where the first variable
has to be thought of the ``position'' variable, whereas the second variable as the ``momentum''
variable.
The Hamiltonian is then
defined as the limit
\be\label{measham}
\caH(\mu, f) =\lim_{N\to\infty, \loc_N(\si)\to\mu} \frac{1}{|\tor|}\left(e^{-|\tor| \lo \bbL_N (\si), f\ra} \loc_N e^{|\tor| \lo \bbL_N (\si), f\ra}\right)
\ee

Note that
$|\tor| \lo \bbL_N (\si), f\ra=\sum_{i\in\tor} \tau_i f(\si)$.

For the computation of \eqref{measham}, we assume $f$ to be a local function. Because the source generator $Q$ is
local we have, that $Q (\tau_k f) =0$ for all $k$ outside the set $D(Q,f)=\{ k: D_Q \cap D_f+k\not=\emptyset\}$.
Therefore, for $\la\subset\Zd$ finite,
\be\label{locprod}
Q \left(\prod_{i\in \la}\tau_i e^f\right) =\left(\prod_{i\in \la\setminus D(Q,f)}\tau_i e^f\right) Q
\left(\prod_{i\in \la\cap D(Q,f)}\tau_i e^ f\right)
\ee
Use \eqref{locprod} to compute
\beq\label{hami}
\caH(\mu, f)&=&  \lim_{N\to\infty, \loc_N(\si)\to\mu} \frac{1}{|\tor|} e^{-\sum_i \tau_i f} \sum_j \tau_j\left(Q e^{\sum_{i} \tau_{i-j} f}\right)
\nonumber\\
&=&
\lim_{N\to\infty, \loc_N(\si)\to\mu} \frac{1}{|\tor|} e^{\sum_i -\tau_i f}\sum_{j\in \tor} \tau_{j}\left(Q e^{\sum_{i\in D(f,Q)+j}\tau_{i-j}f} \right)
e^{\sum_{i\not\in D(f,Q)+j}\tau_i f}
\nonumber\\
&=&
\lim_{N\to\infty, \loc_N(\si)\to\mu} \frac{1}{|\tor|}\sum_{j\in \tor} \tau_j\left( e^{-\sum_{i\in D(f,Q)+j}\tau_{i-j} f }Q e^{\sum_{i\in D(f,Q)+j}\tau_{i-j} f }\right)
\nonumber\\
&=& \lim_{N\to\infty, \loc_N(\si)\to\mu} \frac{1}{|\tor|}
\sum_{j\in \tor} \tau_j\left( e^{-\sum_{k\in D(f,Q)}\tau_{k} f }Q e^{\sum_{k\in D(f,Q)}
\tau_{k} f }\right)
\eeq
We can now introduce the non-linear operator associated to the ``source'' generator $Q$, working on
local functions $f$:
\be\label{kaaq}
\caK_Q f= e^{-\sum_{k\in D(f,Q)} \tau_k f} Q e^{\sum_{k\in D(f,Q)}\tau_k f}
\ee
Using this notation, we obtain from \eqref{hami}
\be\label{hames}
\caH(\mu,f) = \int \caK_Q (f) d\mu
\ee
This Hamiltonian has to be thought of as the analogue of \eqref{hamco} in the present context.

\br
Note that we can write, informally,
\[
\caK_Q f= e^{-\sum_{k\in \Zd} \tau_k f} Q e^{\sum_{k\in\Zd}\tau_k f}
\]
since the terms $k\not\in D(f,Q)$ ``cancel''. This is of course not rigorous
because the infinite sum $\sum_{k\in\Zd}\tau_k f$ does not make sense, but in the ``same way''
as for a formal infinite volume Hamiltonian, where only {\em local energy differences} are well defined.
The advantage of this formal representation is that we clearly see that
$\caK$ is a translation invariant operator, i.e.,
$\caK_Q (f) = \caK_Q (\tau_i f)$, and as a consequence, the Hamiltonian $\caH(\mu, f)$ is translation invariant as well, both in the
measure and in the function, i.e.,
\[
\caH (\tau_k \mu, \tau_r f) = \caH (\mu, f)
\]
for all $k, r\in \Zd$.
Another advantage is that one clearly sees the analogy with the corresponding formula for the
empirical distribution
\eqref{hamco}.
\er

The corresponding Lagrangian is then found by Legendre transformation, i.e.,
\be\label{measlag}
\Xi (\mu, \dot{\mu}) = \sup_{f\in \caC(\Omega)} \left(\int f d\dot{\mu} - \caH(\mu, f) \right)
\ee
where $\dot{\mu}$ denotes a translation invariant signed measure of total mass zero, and
$\mu$ a translation invariant probability measure on $\Omega$.

In general, an explicit expression for $\Xi$ cannot be obtained easily. In the examples below we will compute
$\Xi$ quite explicitly for diffusion processes and show a relative entropy interpretation of $\Xi$ in the context
of interacting particle systems (analogue of finite Markov chains in the previous section) and in the context
of interacting diffusions.

\subsection{Interacting particle systems: the Lagrangian}
We now compute $\caK_Q$ for some of the examples discussed before, starting with interacting particle systems.
The local generator is of the form.
\[
Q f= \sum_{\alpha} r_\alpha (T_\alpha f- f)
\]
where $T_\alpha$ are local transformations, which change coordinates only
in a finite set $D_\alpha$ containing the origin.
This gives
\be\label{intfeng}
\caK_Q f = \sum_\alpha r_\alpha \left( e^{\caD_\alpha (f)} -1\right)
\ee
where the operator $\caD_\alpha$ is defined
by
\[
\caD_\alpha f = \sum_{k\in \Zd}\left(T_\alpha (\tau_k f) -\tau_k f\right)
\]
Notice that the sum is in fact a finite sum since $f$ is local, and the transformation $T_\alpha$ is local as well.
Let us now first zoom in into two familiar examples.
\begin{itemize}
\item[a)] {\bf Independent spin-flip}.
For $E= \{-1,1\}$, and for a single transformation $T\si = \si^0$ (spin-flip), we get
\[
\caD (f) = \sum_{k\in -D_f} (\tau_k f (\si^0)- \tau_k f)
\]
for the special functions $f(\si) = H_A (\si) =\prod_{i\in A} \si_i$ we get
\[
\caD (H_A) = \sum_{k\in -A} -2H_{A+k}
\]
as we found before in \cite{efhr2}.
\item[b)] {\bf Symmetric exclusion process}.
For $E= \{0,1\}$, $d=1$ and $T(\si )= \si^{01}$, where $\si^{01}$ denotes exchange of the values at site $0$ and $1$, i.e.,
$(\si^{01})_j = \si_1 \delta_{j,0} + \si_0 \delta_{j,1} + \si_j (1-\delta_{j,0}-\delta_{j,1})$. We have
\[
\caD (f)(\eta) = \sum_{k: k+D_f\cap \{0,1\} \not= \emptyset} f(\tau_k(\eta^{01}))- f(\eta)
\]
Notice that for $f=\eta_0$ we find only two terms contributing to $\caD(f)$:
\begin{eqnarray*}
\caD(f)&=& \left(((\eta^{01}))_0 - (\eta)_0\right) + \left((\tau_{1}(\eta^{01}))_0 - (\tau_{1}(\eta))_0\right)
\\
&=&
\eta_1-\eta_0 + \eta_0-\eta_1=0
\end{eqnarray*}
which corresponds to the fact that the density of particles is conserved in this process.
\end{itemize}

Returning to the general case now,
the Lagrangian associated with \eqref{intfeng} is
\be\label{intfenglag}
\Xi (\mu,\dot{\mu}) =\sup_{f\in \caC(\Omega)}\left(\int f d\dmu- \int\left(\sum_\alpha r_\alpha \left( e^{\caD_\alpha (f)} -1\right)\right)d\mu\right)
\ee
This expression is reminiscent of
\eqref{finlag} in section 3 (empirical distribution for finite Markov chains). Indeed, a similar relative entropy interpretation of this expression
can be given. We will describe this rather informally, as it is quite analogous to the Girsanov formula computation of the section on finite Markov chains.
First we note that for a translation invariant measure $\mu$, its ``derivative at time zero''
$\loc^*\mu$ is formally given by
\[
(\loc^* \mu )(\si)= \sum_i\sum_\alpha \left(r_\alpha (\tau_i\si) \mu (\tau_i T_\alpha \tau_{-i}\si) - r_\alpha (\tau_i\si) \mu(\tau_i\si)\right)
\]
This object is to be interpreted as working on local functions, i.e., as a distribution.

Suppose now we consider modified rates $\tilde{r}_\alpha (\si)= r_\alpha (\si) e^{f(\si)-f(T_\alpha(\si))}$, and the associated modified local generator
$\tilde{Q}  = \sum_\alpha \tilde{r}_\alpha (T_\alpha-I)$, i.e., the {\em same transformations} are applied but now {\em with new rates}.
Then for a given translation invariant signed measure of total mass zero, we look for those modified rates, i.e., choice of $f$,
such that with the starting measure $\mu$ they produce ``derivative at time zero'' equal to $\dot{\mu}$, i.e.,
\[
\dot{\mu} (\si) = \sum_i\sum_\alpha \left(\tilde{r}_\alpha (\tau_i\si) \mu (\tau_i T_\alpha \tau_{-i}\si) - \tilde{r}_\alpha (\tau_i\si) \mu(\tau_i\si)\right)
\]

The Radon Nikodym derivative of the path space measure of the finite-volume
process (in $\tor$) with rates
$\tilde{r}$ w.r.t.\ the process with rates $r$ is given by the Girsanov formula:
\[
\frac{d\pee^{[0,T],N}{\tilde{r}}}{d\pee^{[0,T],N}{\tilde{r}} }
=
\exp \left( \sum_{i\in \tor}\sum_\alpha\left(\int_0^T\log\frac{\tilde{r}^i_\alpha{\si_s}}{r^i_\alpha(\si_s)}dN^{i,\alpha}_s - \int_0^T\left(\tilde{r}^i_\alpha(\si_s)-r^i_\alpha(\si_s)\right)\right)ds\right)
\]
where $r^i_\alpha$, resp. $\tilde{r}^i_\alpha$ denote the rate to flip
from $\si$ to $\tau_i T_\alpha \tau_{-i} (\si)$, i.e., to apply the transformation $T_\alpha$ around
the lattice site $i$, and $N^{i,\alpha}_t$ the corresponding counting process counting
how many transitions $\si$ to $\tau_i T_\alpha \tau_{-i} (\si)$ have happened in the time
interval $[0,t]$.

We then find, analogously to \eqref{relentlag}
that
the Lagrangian is equal to the limit
\[
\Xi (\mu, \dot{\mu})=\lim_{T\to 0 } \frac1{T}\lim_{N\to\infty} \frac{1}{|\tor|}s(\pee^{[0,T]}_{\tilde{r}, N}|\pee^{[0,T]}_{{r}, N})
\]
i.e., the analogue of \eqref{relentlag}, replacing relative entropy by {\em relative entropy density}.
\subsection{Diffusion processes: the Lagrangian.}
For diffusion processes, let us start with the simplest case of independent Brownian motions in $d=1$.
The general case will be analogous, but the quadratic forms appearing there will be less explicit.
The source generator $Q$ is thus given by
\[
Q f(\si) = \frac12\partial^2_0 f(\si)
\]
where we abbreviated $\partial_0$ to denote the partial derivative w.r.t.\ $\si_0$.
As a consequence, for a local function $f$:
\[
\caK_Q f  = \sum_k Q (\tau_k f) + \frac12\left(\sum_k \partial_0 (\tau_k f)\right)^2
\]
and, reminding that the full generator is the sum of shifts of $Q$, we have
\beq\label{difhammeas}
\caH(\mu, f) &=& \int \caK_Q f d\mu = \int \loc f d\mu + \caJ_\mu (f,f)
\nonumber\\
&=& \lo f, \caL^* \mu \ra + \caJ_\mu (f,f)
\eeq
where
\[
\caJ_\mu (f,f) = \frac12 \int\left(\sum_k \partial_0 (\tau_k f)\right)^2 d\mu
\]
is a $\mu$ dependent quadratic form. This quadratic form is the analogue of \eqref{quad}. Hence,
for the Lagrangian we have
\[
\Xi (\mu, \dot{\mu})= \sup_{f}\left( \lo \dot{\mu}-\loc^*\mu, f\ra -\caJ_\mu (f,f)\right)= \caJ_\mu^* (\dot{\mu}-\loc^*\mu,\dot{\mu}-\loc^*\mu)
\]
where $\caJ_\mu^*$ is a dual quadratic form defined via
\be\label{dualqua}
\caJ_\mu^*(\nu,\nu)= \sup_f \left(\lo\nu,f\ra - \caJ_\mu (f,f)\right)
\ee
for $\nu$ a signed measure of total mass zero.
Notice that this indeed defines a quadratic form
because for $\lambda>0$ (and similarly for $\lambda<0$)
\begin{eqnarray*}
\caJ_\mu^* (\lambda\nu,\lambda\nu) &=& \sup_f (\lambda \lo \nu, f\ra - \caJ_\mu (f,f))
\\
&=&\lambda^2\sup_f (\lo\nu,f/\lambda\ra - \caJ_\mu (f/\lambda,f/\lambda))
\\
&=& \lambda^2 \caJ_\mu^* (\nu,\nu)
\end{eqnarray*}
We see in particular that
$\Xi (\mu, \dot{\mu})$ is zero for a solution of the Kolmogorov forward equation, i.e., if $\dot{\mu} = \loc^* \mu$, which shows also in the present context
that the Markovian evolution
of the distribution $\mu$ is an optimal zero cost trajectory.

Finally, let us turn to the general diffusion case. We split $Q$, the source generator, into a first order part and a second order part:
\[
Q= Q_1+ Q_2
\]
where $Q_2$ contains all second order derivatives (variance part of the diffusion), $Q_1$ all first order derivatives
(drift part).
To $Q_2$ is then associated the quadratic form
\be\label{genquad}
\caJ^Q_\mu (f,f) = \int \left(e^{-\sum_k \tau_k f} Q_2 e^{\sum_k \tau_k f} - Q_2\left(\sum_k \tau_k f\right)\right)d\mu
\ee
Notice that this corresponds to the operator carr\'e du champ associated to $Q_2$, i.e.,
\[
\caJ^Q_\mu (f,f)= \Gamma_2^{Q_2}\left[\sum_k \tau_k f\right]
\]
The Lagrangian is then given by
\be\label{laggendif}
\Xi (\mu,\dot{\mu})= (\caJ^Q_\mu)^* ( \dot{\mu}-\loc^*\mu, \dot{\mu}-\loc^*\mu)
\ee
where $(\caJ^Q_\mu)^*$ is the dual quadratic form of $\caJ^Q$ (as in \eqref{dualqua}).

\end{document}